\documentclass[
11pt]{amsart}
\usepackage{amsthm, amsfonts, amssymb}

\theoremstyle{definition}
\newtheorem{ntn}{Notation}[section]
\newtheorem{dfn}[ntn]{Definition}
\theoremstyle{plain}
\newtheorem{lem}[ntn]{Lemma}

\newtheorem{thm}[ntn]{Theorem}

\theoremstyle{remark}
\newtheorem{rem}[ntn]{Remark}
\newtheorem{exa}[ntn]{Example}

\def\floor[#1]{\lfloor #1 \rfloor }

\newcommand{\z}{\mathbb{Z}}

\newcommand{\lan}{\langle}
\newcommand{\ran}{\rangle}

\newcommand{\sr}{{\rm sr}(R)}

\newcommand{\rn}{R^{2n}}
\newcommand{\GL}{\mathit{GL}}
\newcommand{\Spp}{\mathit{Sp}}

\newcommand{\iur}{\mathcal{IU}(R^{2n})}
\newcommand{\mur}{\mathcal{MU}(R^{2n})}
\newcommand{\hur}{\mathcal{HU}(R^{2n})}
\newcommand{\iu}{\mathcal{IU}}
\newcommand{\hu}{\mathcal{HU}}

\renewcommand{\a}{\mathcal{A}}

\renewcommand{\o}{\mathcal{O}}
\renewcommand{\H}{\tilde{H}}

\renewcommand{\u}{\mathcal{U}}

\newcommand{\s}{\Sigma}
\newcommand{\si}{\sigma}

\newcommand{\del}{\delta}

\newcommand{\arr}{\rightarrow}

\newcommand{\se}{\subseteq}
\newcommand{\bs}{\backslash}
\newcommand{\bcu}{\bigcup}

\newcommand{\mt}{\mapsto}

\newcommand{\ep}{\epsilon}
\newcommand{\lam}{\Lambda}
\newcommand{\lmin}{\Lambda_{{\rm min}}}
\newcommand{\lmax}{\Lambda_{{\rm max}}}
\newcommand{\unit}{\mathit{U}^\epsilon_{2n}(R,\Lambda)}
\newcommand{\eunit}{\mathit{EU}^\epsilon_{2n}(R,\Lambda)}
\newcommand{\OO}{\mathit{O}}
\newcommand{\U}{\mathit{U}}
\newcommand{\usr}{{\rm usr}(R)}
\newcommand{\usra}{{\rm usr}(A)}
\newcommand{\spec}{{\rm Mspec}(R)}
\newcommand{\ur}{{\mathcal{U}'}(R^{2n})}
\newcommand{\urr}{{\mathcal{U}'}}

\DeclareMathAlphabet{\mathds}{U}{dsrom}{m}{n}
\DeclareMathAlphabet{\mathsc}{U}{rsfs}{m}{n}

\begin{document}

\title{Homology stability for Unitary groups}
\author{B. Mirzaii}
\author{W. van der Kallen}

\begin{abstract}
In this paper homology stability for unitary
groups over a ring with finite unitary stable rank is established. Homology
stability of symplectic groups and orthogonal groups appears as a
special case of our results.
\end{abstract}

\maketitle
\maketitle

Our motivation for this work has been to prove homology stability
for unitary groups that is as strong as in the
linear case \cite{kall}. For example we want
homology stability of these groups over any commutative finite dimensional
noetherian ring.

This work is the continuation of our previous work in this direction
\cite{m-vdk}.
After that work, we found that our method, with a little modification,
can also be applied in this more general setting.
Thus our aim is now to prove that homology stabilizes of the unitary
groups over rings with finite unitary stable rank. To do so, as in
\cite{m-vdk}, we prove that the poset of isotropic unimodular sequences is
highly connected. Homology stability for
symplectic groups and orthogonal
groups will appear as a special case.

Our approach is similar to \cite{m-vdk}, especially section 6.
We first compare unimodular
sequences of $\rn$ with a suitable hyperbolic basis of $\rn$, and
we will use that to prove that certain posets are highly connected.
For this we need the higher connectivity of posets of unimodular
sequences due to second author. The higher connectivity of
the poset of isotropic unimodular sequences then follows inductively
as in \cite{m-vdk}.
We conclude with the homology stability theorem.

\section{Posets of unimodular sequences}

Let $R$ be an associative ring with unit. A vector
$(r_1, \dots, r_n) \in R^n$ is called unimodular if there exist
$s_1, \dots, s_n \in R$ such that $\s_{i=1}^ns_ir_i=1$, or
equivalently if the submodule generated by this vector is a free
summand of
the right $R$-module $R^n$. We denote the standard basis of
$R^n$ by $e_1, \dots, e_n$. If $n \le m$, we assume that
$R^n$ is the submodule of $R^m$ generated by $e_1, \dots, e_n \in R^m$.
We denote the set of all unimodular elements of $R^n$ by ${\rm Um(R^n)}$.

We say that a ring $R$ satisfies the {\it stable range condition}
$({\rm S}_m)$,  if $m\geq1$ is an integer so that
for every unimodular vector $(r_1, \dots ,r_m, r_{m+1}) \in R^{m+1}$,
there exist $t_1, \dots ,t_m$ in $R$ such that
$(r_1+t_1r_{m+1}, \dots ,r_m+t_mr_{m+1}) \in R^m$ is unimodular. We say
that $R$ has {\it stable rank} $m$, we denote it with ${\rm sr}(R)=m$,
if $m$ is
the least number such that $({\rm S}_m)$ holds. If such a number does not
exist we say that ${\rm sr}(R)=\infty$.

We refer to \cite{m-vdk} or \cite{kall} for notations such as $\o(V)$.
In particular $\u (R^n)$ denotes the
subposet of $\o(R^n)$ consisting of unimodular sequences. Recall that
a sequence of vectors $v_1, \dots, v_k$ in $R^n$ is
called unimodular when $v_1, \dots, v_k$ is basis of a free direct
summand of $R^n$. Note that if $(v_1, \dots, v_k) \in \o(R^n)$ and if
$n \le m$, it is the
same to say that $(v_1, \dots, v_k)$ is unimodular as a sequence of vectors
in $R^n$ or as a sequence of vectors in $R^m$. We call an element
$(v_1, \dots, v_k)$ of $\u(R^n)$ a $k$-{\it frame}.

\begin{thm}[Van der Kallen]\label{kal5}
Let $R$ be a ring with ${\rm sr}(R) < \infty$ and $n \le m$.
Let $\del$ be $0$ or $1$. Then
\par {\rm (i)} $\o(R^n +  e_{n+1}\del) \cap \u(R^m)$
is $(n-{\rm sr}(R)-1)$-connected.
\par {\rm (ii)} $\o(R^n +  e_{n+1}\del) \cap \u(R^m)_v$ is
$(n-{\rm sr}(R)-|v|-1)$-connected for all $v \in \u(R^m)$.
\end{thm}
\begin{proof}
See \cite[Thm. 2.6]{kall}.
\end{proof}

\section{Hyperbolic spaces and some posets}

{}Let there be an involution on $R$, that is an automorphism of the
additive group of $R$, $R \arr R$ with $r \mt \overline{r}$,
such that $\overline{\overline{r}}=r$ and
$\overline{rs}=\overline{s}\ \overline{r}$. Let $\ep$ be an element in the
center of $R$ such that $\ep\overline{\ep}=1$.
Set $R_\ep =:\{ r-\ep\overline{r} : r \in R\}$ and
$R^\ep =:\{ r \in R: \ep\overline{r}=-r \}$ and observe that
$R_\ep \se R^\ep$. A {\it form parameter}
relative to the involution and $\ep$
is a subgroup $\lam$ of $(R, +)$ such that $R_\ep \se \lam \se R^\ep$ and
$\overline{r}\lam r \se \lam$, for all $r \in R$. Notice that $R_\ep$
and $R^\ep$ are form parameters. We denote them by $\lmin$ and $\lmax$,
respectively. If there is an $s$ in the center of $R$ such that
$s+\overline{s} \in R^\ast$, in particular if $2 \in R^\ast$ , then
$\lmin=\lmax$.

Let $e_{i, j}(r)$ be the $2n \times 2n$-matrix with $r \in R$ in
the $(i, j)$ place and zero elsewhere.
Consider $Q_n=\s_{i=1}^ne_{2i-1, 2i}(1) \in M_{2n}(R)$
and $F_n=Q_n+\ep\ {}^tQ_n =
\s_{i=1}^n(e_{2i-1, 2i}(1)+e_{2i, 2i-1}(\ep)) \in \GL_{2n}(R)$.
Define the
bilinear map $h: R^{2n} \times R^{2n} \arr R$ by
$h(x, y)=
\s_{i=1}^n(\overline{x_{2i-1}}y_{2i}+\ep \overline{x_{2i}}y_{2i-1})$
and $q: \rn \arr R/\lam$ by
$q(x)=
\s_{i=1}^n\overline{x_{2i-1}}x_{2i}\ {\rm mod} \ \lam$, where
$x=(x_1, \dots, x_{2n})$, $y=(y_1, \dots, y_{2n})$
and $\overline{x}=(\overline{x_1}, \dots, \overline{x_{2n}})$.
The triple $(\rn, h,q)$ is called a {\it hyperbolic space}.
By definition the {\it unitary group} relative $\lam$
is the group
\begin{gather*}
\unit:=\{ A \in \GL_{2n}(R) : h(Ax, Ay)=h(x, y), q(Ax)=q(x)
, x, y \in R\}.
\end{gather*}
For more general definitions and the properties of these spaces and groups
see \cite{knus}.

\begin{exa}
(i) Let $\ep=-1$ and let the involution be the identity map ${\rm id}_R$,
then $\lmax=R$. If $\lam=\lmax=R$ then $\unit =:\Spp_{2n}(R)$ is the usual
symplectic group. Note that $R$ is commutative in this case.
\par (ii) Let $\ep=1$ and let the involution be the identity map
${\rm id}_R$, then $\lmin=0$. If $\lam=\lmin=0$ then
$\unit =:\OO_{2n}(R)$ is the usual orthogonal group. As in the
symplectic case,
$R$ is necessarily commutative.
\par (iii) Let $\ep=-1$ and the involution is not the
identity map ${\rm id}_R$.
If $\lam=\lmax$ then $\unit =: \U_{2n}(R)$ is the classical
unitary group corresponding to the involution.
\par (iv) If $\lam =\lmax=R$  then
$\unit=\{ A \in \GL_{2n}(R) : h(Ax, Ay)=h(x, y)\ {\rm for \
all}\ x, y \in R \}= \{ A \in \GL_{2n}(R) : {}^t\overline{A}F_nA=F_n \}$.
\par (v) If $\lam=\lmin=0$ then $\unit=\{ A \in \GL_{2n}(R) :
q(Ax)=q(x)\ {\rm for \ all}\ x\in R \}=
\{ A \in \GL_{2n}(R) : {}^t\overline{A}Q_nA=Q_n\}$.
\end{exa}

Let $\si$ be the permutation of the set of natural numbers given by
$\si(2i)=2i-1$
and  $\si(2i-1)=2i$. For $1 \le i, j \le 2n$, $i \neq j$, and every
$r \in R$ define
\begin{gather*}
E_{i, j}(r)=
\begin{cases}
I_{2n}+ e_{i, j}(r) & \text{if $i=2k-1, j=\si(i), r \in \lam$}\\
I_{2n}+ e_{i, j}(r) & \text{if $i=2k, j=\si(i), \overline{r} \in \lam$}\\
I_{2n}+ e_{i, j}(r)+ e_{\si(j), \si(i)}(-\overline{r})&
\text{if $i+j=2k$, $i\neq j$}\\
I_{2n}+ e_{i, j}(r)+ e_{\si(j), \si(i)}(-\ep^{-1}\overline{r})&
\text{if $i \neq\si(j)$, $i=2k-1$, $j=2l$}\\
I_{2n}+ e_{i, j}(r)+ e_{\si(j), \si(i)}(\ep\overline{r})&
\text{if $i \neq\si(j)$, $i=2k$, $j=2l-1$} \end{cases}
\end{gather*}
where $I_{2n}$ is the identity element of $\GL_{2n}(R)$.
It is easy to see that $E_{i, j}(r) \in \unit$. 
Let  $\eunit$ be the group generated by the $E_{i, j}(r)$, $r \in R$.
We call it  {\it elementary unitary group}.

A nonzero vector $x \in \rn$ is called isotropic if $q(x)=0$.
This shows automaticly that if $x$ is isotropic then $h(x, x)=0$.
We say that a subset $S$ of $R^{2n}$ is isotropic
if for every $x \in S$, $q(x)=0$ and for every $x, y \in S$,
$h(x, y)=0$.
If $h(x, y)=0$, then we say that $x$ is perpendicular to $y$.
We denote by $\lan S \ran$ the submodule of $\rn$ generated by $S$,
and by $\lan S \ran^\perp$ the submodule consisting of all the
elements of $\rn$ which are perpendicular to all the elements of $S$.

 From now, we fix an involution, an $\ep$, a form parameter $\lam$
 and we consider the triple $(\rn, h, q)$ as defined above.

\begin{dfn}[Transitivity condition]
Let $r \in R$ and define
$C_r^\ep(\rn, \lam)=\{ x \in {\rm Um}(\rn): q(x)= r\ {\rm mod}\ \lam \}$.
We say that $R$ satisfies the transitivity
condition $({\rm T}_n)$, if $\eunit$
acts transitively on $C_r^\ep(\rn, \lam)$, for every $r \in R$. It is easy
to see that $e_1 + re_2 \in C_r^\ep(\rn, \lam)$.
\end{dfn}

\begin{dfn}[Unitary stable range]
We say that a ring $R$ satisfies
the unitary stable range condition 
if $R$ satisfies the
conditions $({\rm S}_m)$ and $({\rm T}_{m+1})$. We say that $R$ has
unitary stable rank $m$, we denote it with $\usr$, if $m$ is the least
number such that $({\rm S}_m)$ and $({\rm T}_{m+1})$ are satisfied.
If such a number does not exist we say that $\usr=\infty$. Clearly
$\sr \le \usr$.
\end{dfn}

\begin{rem}
Our definition of unitary stable range is a little different than the
one in
\cite[Chap. VI 4.6]{knus}. In fact $\usr +1=m+1={\rm USR}_m$ where
${\rm USR}_m$ is the unitary stable rank as defined in \cite{knus}.
\end{rem}

\begin{exa}
Let $R$ be a commutative noetherian ring where the dimension $d$
of the maximal spectrum
$\spec$ is finite. If $A$ is a finite $R$-algebra then $\usra \le d+2$
(see \cite[Thm 2.8]{vas3}, \cite[Thm. 6.1.4]{knus}).
\end{exa}

\begin{lem}\label{b-w0}
Let $R$ be a ring with $\usr < \infty$. Assume $n \ge \usr + k+1$ and
$(v_1,\dots, v_k) \in \u(\rn)$. Then there is a hyperbolic basis
$\{x_1, y_1, \dots, x_n, y_n\}$ of $\rn$ such that
$v_1, \dots, v_k \in \lan x_1, y_1, \dots, x_k, y_k\ran$.
\end{lem}
\begin{proof}
The proof is by induction on $k$.
If $k=1$, by definition of unitary stable range there is an $E \in \eunit$
such that $Ev_1=e_1 + re_2$. So the base of induction is true.
Let $k \ge 2$ and assume the induction hypothesis.
Arguing as in the base of the induction
 we can assume that $v_1=(1, r, 0, \dots, 0)$. Let
$W=e_2 + \s_{i=2}^ne_{2i}R$. By theorem \ref{kal5}, the poset
$F:=\o(W) \cap \u(\rn)_{(v_1,\dots, v_k)}$ is $((n-1)-\sr -k-1)$-connected.
Since $n \ge \usr + k+1 \ge \sr +k+1$, it follows that $F$ is not empty.
Choose $(w, v_1,\dots, v_k) \in \u(\rn)$ where $w \in W$. Then
$(w, v_1-wr, v_2, \dots, v_k) \in \u(\rn)$. But $(w, v_1-rw)$ is a
hyperbolic pair, so there is an $E \in \eunit$ such that
$Ew=e_{2n-1}, E(v_1- wr)=e_{2n}$ by \cite[Chap. VI, Thm. 4.7.1]{knus}. Let
$(Ew, E(v_1-wr), Ev_2, \dots, Ev_k)=:( w_0, w_1, \dots, w_k)$ where
$w_i=(r_{i, 1}, \dots ,  r_{i, 2n})$. Put
$u_i=w_i-e_{2n-1}r_{i, 2n-1}-e_{2n}r_{i, 2n}$ for $2 \le i \le k$.
Then $(u_2, \dots, u_k) \in \u(R^{2n-2})$. Now by
induction there is a hyperbolic
basis $\{a_2, b_2, \dots, a_n, b_n\}$ of $R^{2n-2}$ such that
$u_i \in  \lan a_2, b_2, \dots, a_k, b_k \ran$. Let $a_1=e_{2n-1}$ and
$b_1=e_{2n}$. Then $w_i \in \lan a_1, b_1, \dots, a_k, b_k \ran$.
But $Ev_1=w_1+Ewr=e_{2n}+e_{2n-1}r$,
$Ev_i=w_i$ for $2 \le i \le k$ and considering
$x_i=E^{-1}a_i$, $y_i=E^{-1}b_i$, one sees that
$v_1, \dots, v_k \in \lan x_1, y_1, \dots, x_k, y_k\ran$.
\end{proof}

Let $Z_n=\{x \in \rn : q(x)=0\}$ and put $\ur=\o(Z_n) \cap \u(\rn)$.

\begin{lem}\label{uu}
Let $R$ be a ring with $\sr < \infty$ and $n \le m$. Then
\par {\rm (i)} $\o(\rn) \cap \urr(R^{2m})$ is $(n-\sr-1)$-connected,
\par {\rm (ii)} $\o(\rn) \cap \urr(R^{2m})_v$ is $(n-\sr-|v|-1)$-connected
for every $v \in \urr(R^{2m})$,
\par {\rm (iii)} $\o(\rn) \cap \urr(R^{2m}) \cap \u(R^{2m})_v$
is $(n-\sr-|v|-1)$-connected
for every $v \in \u(R^{2m})$.
\end{lem}
\begin{proof}
Let $W=\lan e_2, e_4, \dots, e_{2n}\ran$ and $F:=\o(\rn) \cap \urr(R^{2m})$.
It is easy to see that $\o(W) \cap F=\o(W) \cap \u(R^{2m})$ and
$\o(W) \cap F_u=\o(W) \cap \u(R^{2m})_u$ for every $u \in \urr(R^{2m})$.
By theorem \ref{kal5}, the poset $\o(W) \cap F$ is $(n-\sr-1)$-connected and
the poset $\o(W) \cap F_u$ is $(n-\sr-|u|-1)$-connected for every
$u \in F$. It follows from
lemma \cite[2.13 (i)]{kall} that $F$ is
$(n-\sr-1)$-connected. The proof of (ii) and (iii)
is similar to the proof of (i).
\end{proof}

\begin{lem}\label{u-i}
Let $R$ be a ring with $\usr < \infty$ and let  $(v_1, \dots, v_k) \in \ur$.
If $n \ge \usr + k+1$ then
$\o(\langle v_1, \dots, v_k\rangle^\perp) \cap \ur_{(v_1, \dots, v_k)}$
is $(n-\usr-k-1)$-connected.
\end{lem}
\begin{proof}
By lemma \ref{b-w0} there is a hyperbolic basis
$\{x_1, y_1, \dots, x_n, y_n\}$ of $\rn$ such that
$v_1, \dots, v_k \in \lan x_1, y_1, \dots, x_k, y_k\ran$. Let
$W=\lan x_{k+1}, y_{k+1}, \dots, x_n, y_n\ran \simeq R^{2(n-k)}$ and
$F:=\o(\lan v_1, \dots, v_k\ran^\perp) \cap \ur_{(v_1, \dots, v_k)}$.
It is easy to see that $\o(W) \cap F=\o(W) \cap \ur$.
Let $V=\lan v_1, \dots, v_k \ran$, then
$\lan x_1, y_1, \dots, x_k, y_k\ran = V \oplus P$ where P is a
(finitely generated) projective module. Consider
$(u_1, \dots, u_l) \in F \bs \o(W)$ and let $u_i=x_i+y_i$ where
$x_i \in V$ and $y_i \in P \oplus W$. One should
notice that $(u_1-x_1, \dots, u_l-x_l) \in \u(\rn)$
and not necessarily in $\ur$. It is not difficult to see that
$\o(W) \cap F_{(u_1, \dots, u_l)}=\o(W)
\cap \ur \cap \u(\rn)_{(u_1-x_1, \dots, u_l-x_l)}$.
By lemma \ref{uu}, $\o(W) \cap F$ is $(n-k-\usr-1)$-connected
and $\o(W) \cap F_u$ is $(n-k-\usr-|u|-1)$-connected for every
$u \in F\bs \o(W)$. It follows
from lemma \cite[2.13 (i)]{kall} that $F$ is
$(n-\usr-k-1)$-connected.
\end{proof}

\section{Posets of isotropic and hyperbolic unimodular sequences}

Let $\iur$ be the set of sequences $(x_1, \dots, x_k)$, $x_i \in R^{2n}$,
such that $x_1, \dots, x_k$ form a basis for an isotropic direct summand
of $R^{2n}$. Let
$\hur$ be the set of sequences $((x_1, y_1), \dots, (x_k, y_k))$
such that $(x_1, \dots, x_k)$, $(y_1, \dots, y_k) \in \iur$,
$h(x_i, y_j)=\delta_{i, j}$, where $\delta_{i, j}$ is the Kronecker delta.
We call $\iur$ and $\hur$ the
poset of isotropic unimodular sequences and the poset of
hyperbolic unimodular sequences, respectively.
For $1 \le k \le n$, let
$\iu(\rn, k)$ and $\hu(\rn, k)$ be the set of all elements of length $k$
of $\iur$ and $\hur$ respectively.
We call the elements of $\iu(\rn, k)$ and $\hu(\rn, k)$ the
isotropic $k$-frames and the hyperbolic $k$-frames, respectively.
Define the poset $\mur$ as the set of
$((x_1, y_1), \dots, (x_k, y_k)) \in \o(\rn \times \rn)$ such that,
{\rm (i)} $(x_1, \dots, x_k) \in \iur$,
{\rm (ii)} for each $i$, either $y_i=0$ or $(x_j, y_i)=\delta_{ji}$,
{\rm (iii)} $\lan y_1, \dots, y_k \ran$ is isotropic.
We identify
$\iur$ with $\mur \cap \o(\rn \times \{0\})$ and
$\hur$ with $\mur \cap \o(\rn \times (\rn \bs \{0\}))$.

\begin{lem}\label{vas0}
Let $R$ be a ring
with $\usr < \infty$. If $n \ge \usr + k$ then
$\eunit$ acts transitively on $\iu(\rn, k)$ and $\hu(\rn, k)$.
\end{lem}
\begin{proof}
The proof is by induction on $k$. If $k=1$,
by definition $\eunit$ acts transitively on $\iu(\rn, 1)$ and by
\cite[Chap. VI, Thm. 4.7.1]{knus} the group $\eunit$ acts transitively on
$\hu(\rn, 1)$.
The rest is an easy induction and the fact that for every isotropic $k$-frame
$(x_1, \dots, x_k)$ there is an  isotropic $k$-frame
$(y_1, \dots, y_k)$ such that $((x_1, y_1), \dots, (x_k, y_k))$ is a
hyperbolic $k$-frame \cite[Chap. I, Cor. 3.7.4]{knus}.
\end{proof}

\begin{lem}\label{charn}
Let $R$ be a ring with $\usr < \infty$, and let $n \ge \usr + k$.
Let $((x_1, y_1), \dots, (x_k, y_k))$ $\in \hur$, $(x_1, \dots, x_k)$ $ \in
\iur$
and $V=\lan x_1, \dots, x_k\ran$. Then
\par {\rm(i)} $\iur_{(x_1, \dots, x_k)} \simeq \iu(R^{2(n-k)})\lan V \ran$,
\par {\rm(ii)} $\hur \cap \mur_{((x_1, 0), \dots, (x_k, 0))}
\simeq \hur_{((x_1, y_1), \dots, (x_k, y_k))}\lan V \times V \ran$,
\par {\rm(iii)} $\hur_{((x_1, y_1), \dots, (x_k, y_k))} \simeq
\hu(R^{2(n-k)})$.
\end{lem}
\begin{proof}
See \cite{ch}, the proof of lemma 3.4 and the proof of Thm. 3.2.
\end{proof}

For a real number $l$, by $\floor [l]$ we mean the largest integer $n$ with
$n\leq l$.

\begin{thm}\label{b-w1}
The poset $\iur$ is
$\floor [\frac{n-\usr-3}{2}]$-connected and $\iur_x$ is
$\floor [\frac{n-\usr-|x|-3}{2}]$-connected for every
$x \in \iur$.
\end{thm}
\begin{proof}
If $n\leq \usr+1$, the result is clear,
so let $n>\usr+1$. Let
$X_v=\iur \cap \ur_v \cap \o(\langle v \rangle^\perp)$,
for every $v \in \ur$, and put $X:= \bcu_{v \in F}X_v$ where $F=\ur$.
It follows from lemma \ref{vas0}
that $\iur_{\le n-\usr-1} \se  X$.
So to treat $\iur$, it is enough to prove that
$X$ is $\floor [\frac{n-\usr-3}{2}]$-connected.
First we prove that  $X_v$ is
$\floor [\frac{n-\usr-|v|-3}{2}]$-connected for every $v \in F$.
The proof is by descending
induction on $|v|$.
If $|v| > n-\usr-1$, then $\floor [\frac{n-\usr-|v|-3}{2}]<-1$.
In this case there
is nothing to prove. If $n-\usr-2\le|v| \le n-\usr-1$, then
$\floor [\frac{n-\usr-|v|-3}{2}]= -1$,
so we must prove that $X_v$ is nonempty. This follows from lemma \ref{b-w0}.
Now assume $|v|\leq n-\usr-3$ and assume
by induction that $X_w$ is $\floor [\frac{n-\usr-|w|-3}{2}]$-connected
for  every $w$, with $|w| >|v|$.
Let $l=\floor [\frac{n-\usr-|v|-3}{2}]$, and observe that
$ n-|v|-\usr-1\geq l+2$.
Put $T_w=\iur \cap \ur_{wv} \cap \o(\lan wv \ran^\perp)$
where $w \in G_v=\ur_v\cap \o(\langle v \rangle^\perp)$
and put $T:=\bcu_{w \in G_v}T_w$ . It follows by lemma \ref{b-w0}
that $(X_v)_{\le n-|v|-\usr-1} \se T$.
So it is enough to prove that $T$ is $l$-connected.
The poset $G_v$ is $l$-connected by lemma \ref{u-i}.
By induction,
$T_w$ is $\floor [\frac{n-\usr-|v|-|w|-3}{2}]$-connected.
But ${\rm min}\{l-1, l-|w|+1\} \le \floor [\frac{n-\usr-|v|-|w|-3}{2}]$,
so $T_w$ is
${\rm min}\{l-1, l-|w|+1\}$-connected. For every $ y \in T$,
$\a_y= \{w \in G_v: y \in T_w\}$ is isomorphic to
$\ur_{vy} \cap \o(\langle
vy \rangle^\perp)$ so by lemma \ref{u-i},
it is $(l-|y|+1)$-connected.
Let $w\in G_v$ with $|w|=1$.
For every $z\in T_w$ we have $wz\in X_v$,
so $T_w$ is contained in a cone,  call it
$C_w$, inside $X_v$. Put $C(T_w)=T_w \cup (C_w)_{\le n-|v|-\usr-1}$.
Thus $C(T_w) \se T$. The poset $C(T_w)$ is $l$-connected because
$C(T_w)_{\le n-|v|-\usr-1} = (C_w)_{\le n-|v|-\usr-1}$.
Now by theorems \ref{kal5} and \cite[4.7]{m-vdk}
, $T$ is $l$-connected.
In other words, we have now shown that
$X_v$ is $\floor [\frac{n-\usr-|v|-3}{2}]$-connected.
By knowing this one can prove, in a similar way, that $X$ is
$\floor [\frac{n-\usr-3}{2}]$-connected. (Just pretend that $|v|=0$.)

Now consider the poset $\iur_x$ for an $x=(x_1, \dots, x_k) \in \iur$.
The proof is by induction on $n$. If $n=1$, everything is easy.
Similarly, we may assume $n-\usr-|x|-1\geq 0$.
Let $l=\floor [\frac{n-\usr-|x|-3}{2}]$. By lemma \ref{charn},
$\iur_x \simeq \iu(R^{2(n-|x|)})\lan V \ran$, where
$V=\lan x_1, \dots, x_k \ran$. In the above we proved that $\iu(R^{2(n-|x|)})$
is $l$-connected and by induction, the poset
$\iu(R^{2(n-|x|)})_y$ is $\floor [\frac{n-|x|-\usr-|y|-3}{2}]$-connected
for every $y \in \iu(R^{2(n-|x|)})$.
But $l-|y| \le \floor [\frac{n-|x|-\usr-|y|-3}{2}]$. So
$\iu(R^{2(n-|x|)})\lan V \ran$ is $l$-connected  by
lemma \cite[4.1]{m-vdk}
. Therefore $\iur_x$
is $l$-connected.
\end{proof}

\begin{thm}\label{b-w2}
The poset $\hur$ is $\floor [\frac{n-\usr-4}{2}]$-connected and
$\hur_x$ is $\floor [\frac{n-\usr-|x|-4}{2}]$-connected for every
$x \in \hur$.
\end{thm}
\begin{proof}
The proof is by induction on $n$. If $n=1$, then everything is trivial.
Let $F=\iur$ and $X_v=\hur \cap \mur_v$,
for every $v \in F$. Put  $X:=\bcu_{v \in F} X_v$.
It follows from lemma \ref{vas0} that $\hur_{\le n-\usr-1} \se X$.
Thus to treat $\hur$, it is enough to prove that $X$ is
$\floor [\frac{n-\usr-4}{2}]$-connected, and we may assume $n\geq\usr+2$.
Take
$l=\floor [\frac{n-\usr-4}{2}]$ and $V=\lan v_1, \dots, v_k \ran$, where
$v=(v_1, \dots, v_k)$. By lemma \ref{charn}, there is an isomorphism
$X_v \simeq \hu(R^{2(n-|v|)}) \lan V \times V \ran$, if $n\geq\sr+|v|$.
By induction
$\hu(R^{2(n-|v|)})$ is $\floor [\frac{n-|v|-\usr-4}{2}]$-connected and
again by induction
$\hu(R^{2(n-|v|)})_y$ is
$\floor [\frac{n-|v|-\usr-|y|-4}{2}]$-connected for every
$y \in \hu(R^{2(n-|v|)})$. So by  lemma \cite[4.1]{m-vdk}
, $X_v$ is
$\floor [\frac{n-|v|-\usr-4}{2}]$-connected.
Thus the poset $X_v$ is
${\rm min}\{l-1, l-|v|+1\}$-connected.
Let
$x=((x_1, y_1), \dots, (x_k, y_k))$. It is easy to see that
$\a_x=\{ v \in F: x \in X_v\} \simeq \iur_{(x_1, \dots, x_k)}$.
By the above theorem \ref{b-w1}, $\a_x$ is
$\floor [\frac{n-\usr-k-3}{2}]$-connected. But
$l-|x|+1 \le \floor [\frac{n-\usr-k-3}{2}]$,
so $\a_x$ is $(l-|x|+1)$-connected.
Let $v=(v_1) \in F$, $|v|=1$, and let
$D_v:= \hur_{(v_1, w_1)}\simeq \hu(R^{2(n-1)})$ where $w_1 \in \rn$ is a
hyperbolic
dual of $v_1 \in \rn$.
Then $D_v \se X_v$ and $D_v$ is contained in a
cone, call it $C_v$, inside $\hur$. Take
$C(D_v):=D_v \cup (C_v)_{\le n-\usr-1}$.
By induction $D_v$ is $\floor [\frac{n-1-\usr-4}{2}]$-connected and so
$(l-1)$-connected.
Let $Y_v= X_v \cup C(D_v)$. By the Mayer-Vietoris
theorem and the fact that $C(D_v)$ is $l$-connected,
we get the exact sequence
\[
\H_l(D_v, \z) \overset{(i_v)_\ast}{\arr} \H_l(X_v, \z)
\arr \H_l(Y_v, \z) \arr 0.
\]
where $i_v: D_v \arr X_v$ is the inclusion.
By induction $(D_v)_w$ is $\floor [\frac{n-1-\usr-|w|-4}{2}]$-connected
and so
$(l-|w|)$-connected, for $w\in D_v$.
By lemma \cite[4.1(i)]{m-vdk} 
and lemma \ref{charn},
$(i_v)_\ast$ is an isomorphism,
and by exactness of the above sequence we get
$\H_l(Y_v, \z)=0$. If $l \ge 1$ by the Van Kampen theorem
$\pi_1(Y_v, x) \simeq \pi_1(X_v, x)/N$ where $x \in D_v$ and $N$ is the
normal subgroup generated by the image of the map
$(i_v)_\ast: \pi_1(D_v, x) \arr \pi_1(X_v, x)$. Now by lemma
\cite[4.1(ii)]{m-vdk}
, $\pi_1(Y_v, x)$ is trivial.
Thus by the Hurewicz theorem \cite[2.1]{m-vdk}
, $Y_v$ is $l$-connected.
By having all this we can apply  theorem
\cite[4.7]{m-vdk}
and so $X$ is $l$-connected.
The fact that $\hur_x$ is
$\floor [\frac{n-\usr-|x|-4}{2}]$-connected follows from the above
and lemma \ref{charn}.
\end{proof}

\section{ Homology stability}

{}From theorem \ref{b-w2} one can get the homology stability of
unitary groups as Charney proved in \cite[Sec. 4]{ch}. Here we only
formulate the theorem and for the proof we refer to Charney's paper.

\begin{thm}
Let $R$ be a ring with $ \usr < \infty$. Then for every abelian
group $L$ the homomorphism
$ {\psi_n}_\ast : H_i(\unit, L) \arr H_i(U^\ep_{2n+2}(R, \lam), L)$
is surjective for $n \ge 2i+\usr+3$ and bijective for
$n \ge 2i+\usr+4$, where
$ {\psi_n} : \unit \arr {U}^\ep_{2n+2}(R, \lam)$, $A \mt
\left(\begin{array}{cc}
A & 0      \\
0 &  I_2
\end{array} \right)$.
\end{thm}
\begin{proof}
See \cite[Sec. 4]{ch}.
\end{proof}

\

 e-mail:\quad \texttt{mirzaii@math.uu.nl\qquad vdkallen@math.uu.nl}

\begin{thebibliography}{99}








\bibitem{ch} Charney, R. A generalization of a theorem of Vogtmann.
J. Pure Appl. Algebra {\bf 44} (1987), 107--125.






\bibitem{knus} Knus, M. A. Quadratic and Hermitian forms over rings.
Grundlehren der Mathematischen Wissenschaften , 294. Springer-Verlag,
Berlin, 1991.




\bibitem{m-vdk} Mirzaii, B.; Van der Kallen, W. Homology stability
for symplectic groups. 
{\tt arXiv:math.KT/0110163}











\bibitem{kall} Van der Kallen, W. Homology stability for linear groups.
Invent. Math. {\bf 60} (1980), 269--295.

\bibitem{vas3} Vaserstein, L. N. Stabilization of unitary and orthogonal
groups over a ring with involution.
Math. USSR Sbornik {\bf 10} (1970),
no.~3, 307--326.





\end{thebibliography}
\end{document}